\documentclass[11pt]{amsart}

\usepackage{amsfonts,amssymb,stmaryrd,amscd,amsmath,latexsym,amsbsy}

\usepackage{amssymb}
\usepackage{amsfonts}
\usepackage{latexsym}

\newtheorem{theorem}{Theorem}[section]

\newtheorem{proposition}[theorem]{Proposition}
\newtheorem{corollary}[theorem]{Corollary}
\theoremstyle{definition}
\newtheorem{definition}[theorem]{Definition}

\newtheorem{example}[theorem]{Example}

\newtheorem{conjecture}[theorem]{Conjecture}

\newcommand{\h}{\mathfrak{h}}

\newcommand{\ben}{\begin{enumerate}}
\newcommand{\een}{\end{enumerate}}

\theoremstyle{plain}

\newtheorem*{sol}{Solution}

\theoremstyle{definition}

\theoremstyle{remark}

\newcommand{\solu}[1]{\begin{sol}{\bf (\ref{#1})}}

\begin{document}

\centerline{\Large\bf Exploring noncommutative algebras 
via deformation theory}

\vskip .1in

\centerline{\Large\bf Pavel Etingof (MIT)}

\vskip .1in

In this lecture \footnote{This lecture was delivered 
at ``Giornata IndAM'', Naples, June 7, 2005. I would like to thank 
the organizers, in particular Prof. Corrado 
De Concini and Paolo Piazza for this wonderful opportunity.
I am also grateful to J. Stasheff for useful comments.}
I would like to address the following question: 
given an associative algebra $A_0$, what are the possible ways to deform
it? Consideration of this question for concrete algebras often leads 
to interesting mathematical discoveries. I will discuss several approaches
to this question, and examples of applying them.  

\section{Deformation theory} 

\subsection{Formal deformations} 

The most general approach to the question ``how to deform $A_0$?'' 
is the theory of formal deformations. 

Let $k$ be a field and $K:=k[[\hbar_1,...,\hbar_\ell]]$ the ring 
of formal power series in variables $\hbar_i$. Let ${\frak m}$ be the maximal ideal in $K$. 

A $K$-module $M$ is 
said to be {\bf topologically free} if it is isomorphic 
to $M_0[[\hbar_1,...,\hbar_\ell]]$ for some vector space $M_0$. 

Let $A_0$ be an algebra over $k$.\footnote{By ``an algebra'' we always mean 
an associative algebra with unit.} 

\begin{definition} An {\bf $\ell$-parameter 
flat formal deformation} of $A_0$ is 
an algebra $A$ over $K$ which is topologically free as a $K$-module,
together with an isomorphism of algebras 
$\phi: A/{\frak m}\to A_0$.\footnote{The word ``flat'' refers to the fact that 
$A$ is a (topologically) flat module over $K$, i.e. the functor of completed tensor product 
with this module is exact.}  
\end{definition}

For simplicity we will mostly consider 1-parameter deformations. 
If $A$ is a 1-parameter flat formal deformation of $A_0$ then 
we can choose an identification $A\to A_0[[\hbar]]$ 
as $K$-modules, which reduces to $\phi$ modulo $\hbar$. 
Then the algebra structure on $A$ transforms into a new 
$K$-linear multiplication law $\mu$ on $A_0[[\hbar]]$. 
Such a multiplication law is determined by 
the product $\mu(a,b)$, $a,b\in A_0\subset A_0[[\hbar]]$, 
which is given by the formula 
$$
\mu(a,b)=\mu_0(a,b)+\hbar \mu_1(a,b)+\hbar^2\mu_2(a,b)+..., a,b\in A_0,
$$ 
where $\mu_i: A_0\otimes A_0\to A_0$ are linear maps, and 
$\mu_0(a,b)$ is the undeformed product $ab$ in $A_0$. 
Thus, to find formal deformations 
of $A_0$ means to find all such series $\mu$ which satisfy the associativity 
equation, modulo the automorphisms of the $K$-module $A_0[[\hbar]]$
which are the identity modulo $\hbar$. 
\footnote{Note that we don't have to worry about the existence of a unit in 
$A$ since a flat formal deformation of an algebra with unit always has a 
unit}  
 
The associativity equation 
$\mu\circ(\mu\otimes Id)=\mu\circ (Id\otimes \mu)$ 
reduces to a hierarchy of linear equations: 
\begin{equation}\label{asso}
\sum_{s=0}^N \mu_s(\mu_{N-s}(a,b),c)=
\sum_{s=0}^N \mu_s(a,\mu_{N-s}(b,c)).
\end{equation}
(These equations are linear in $\mu_N$ if $\mu_i$, $i<N$, are known).

\subsection{Hochschild cohomology}

These equations can be analyzed using Hochschild cohomology. 
Let us recall its definition. Let $M$ be a bimodule over $A_0$.  
A Hochschild $n$-cochain of $A_0$ with coefficients in $M$ 
is a linear map $A_0^{\otimes n}\to M$. 
The space of such cochains is denoted by $C^n(A_0,M)$. 
The differential $d:C^n(A_0,M)\to C^{n+1}(A_0,M)$ is defined 
by the formula 
$$
df(a_1,...,a_{n+1})=f(a_1,...,a_n)a_{n+1}-f(a_1,...,a_na_{n+1})
$$
$$
+f(a_1,a_2a_3,...,a_{n+1})-...+(-1)^nf(a_1a_2,...,a_{n+1})+
(-1)^{n+1}a_1f(a_2,...,a_{n+1}). 
$$
It is easy to show that $d^2=0$, and one defines the Hochschild cohomology 
$H^\bullet(A_0,M)$ to be the cohomology of the complex 
$(C^\bullet(A_0,M),d)$. If $M=A_0$, the algebra itself, then we will 
denote $H^\bullet(A_0,M)$ by $H^\bullet(A_0)$ (it is an algebra). 
For example, $H^0(A_0)$ is the center of $A_0$, and 
$H^1(A_0)$ is the quotient of the Lie algebra of derivations of $A_0$ by inner derivations.  

The following are standard facts from deformation theory
(due to Gerstenhaber \cite{Ge}), which can be
checked directly.  

1.  The linear equation for $\mu_1$ says that $\mu_1$ is a Hochschild 
2-cocycle. Thus algebra structures on $A_0[\hbar]/\hbar^2$ 
deforming $\mu_0$ are parametrized by the space $Z^2(A_0)$ of 
Hochschild 2-cocycles of $A_0$ with values in $M=A_0$.  

2. If $\mu_1,\mu_1'$ are two 2-cocycles such that 
$\mu_1-\mu_1'$ is a coboundary, then the algebra structures 
on $A_0[\hbar]/\hbar^2$ corresponding to $\mu_1$ and $\mu_1'$ are 
equivalent by a transformation of $A_0[\hbar]/\hbar^2$ 
that equals the identity 
modulo $\hbar$, and vice versa. Thus equivalence classes 
of multiplications on $A_0[\hbar]/\hbar^2$ deforming $\mu_0$
are parametrized by the cohomology $H^2(A_0)$. 

3. The linear equation for $\mu_N$ says that $d\mu_N$ 
is a certain quadratic expression $b_N$ in $\mu_0,\mu_1,...,\mu_{N-1}$. 
This expression is always a Hochschild 3-cocycle, and the equation is 
solvable iff it is a coboundary. Thus the cohomology class 
of $b_N$ in $H^3(A_0)$ is the only obstruction to solving this equation. 

\subsection{Universal deformation}

In particular, if $H^3(A_0)=0$ then the equation for $\mu_n$ can be solved 
for all $n$, and for each $n$ the freedom in choosing the solution, 
modulo equivalences, is the space $H:=H^2(A_0)$. Thus there exists 
an algebra structure over $k[[H]]$ on the space $A_u:=A_0[[H]]$ 
of formal functions from $H$ to $A_0$, $a,b\mapsto \mu_u(a,b)\in A_0[[H]]$,
($a,b\in A_0$), such that $\mu_u(a,b)(0)=ab\in A_0$, and 
every 1-parameter flat formal deformation $A$ of $A_0$ is given 
by the formula $\mu(a,b)(\hbar)=\mu_u(a,b)(\gamma(\hbar))$
for a unique formal series $\gamma\in \hbar H[[\hbar]]$,
with the property that $\gamma'(0)$ is the cohomology class 
of the cocycle $\mu_1$. 

Such an algebra $A_u$ is called a universal deformation of $A_0$. 
It is unique up to an isomorphism. 

Thus in the case $H^3(A_0)=0$, deformation theory allows us 
to completely classify 1-parameter flat formal deformations of $A_0$.
In particular, we see that the ``moduli space'' parametrizing formal deformations of $A_0$
is a smooth space -- it is the formal neighborhood of zero in $H$.

\subsection{Quantization of Poisson structures}

If $H^3(A_0)$ is nonzero then in general the universal 
deformation parametrized by $H$ does not exist, as there are obstructions
to deformations. In this case, the moduli space of deformations 
will be a closed subscheme of $H$, which is often singular. 
On the other hand, even when $H^3(A_0)\ne 0$, the universal 
deformation parametrized by $H$ may exist (although 
it may be more difficult to prove than in the vanishing case).
In this case one says that the deformations of $A_0$ are {\bf unobstructed}
(since all obstructions vanish even though the space of obstructions doesn't).

To illustrate these statements, consider the quantization theory of Poisson manifolds. 
Let $M$ be a smooth $C^\infty$-manifold or a smooth affine algebraic variety over $\Bbb C$, 
and $A_0$ the structure algebra of $M$.

{\bf Remark.} In the $C^\infty$-case, we will consider only local maps 
$A_0^{\otimes n}\to A_0$, i.e. those given by polydifferential operators, and 
all deformations and the Hochschild cohomology 
is defined using local, rather than general, cochains. 

\begin{theorem}(Hochschild-Kostant-Rosenberg) \cite{HKR} 
$H^i(A_0)=\Gamma(M,\wedge^iTM)$ as a module over $A_0=H^0(A_0)$. 
\end{theorem}

In particular, $H^2$ is the space of bivector fields, and $H^3$ the space of trivector fields. 
So the cohomology class of $\mu_1$ is a bivector field; in fact, it is $\pi(a,b):=\mu_1(a,b)-\mu_1(b,a)$, 
since any 2-coboundary in this case is symmetric. The equation for $\mu_2$ says that $d\mu_2$
is a certain trivector field that depends quadratically on $\pi$. It is easy to show that this is 
the Schouten bracket $[\pi,\pi]$. Thus, for the existence of $\mu_2$ it is necessary that $[\pi,\pi]=0$,
i.e. that $\pi$ be a {\bf Poisson bracket}.

Suppose now that $\pi$ is a Poisson bracket, i.e. $[\pi,\pi]=0$. In this 
case the algebra $A=A_0[[\hbar]]$ with the 
product $\mu$ is said to be {\bf a quantization} of $\pi$, and $(M,\pi)$ the {\bf quasiclassical limit} of 
$(A,\mu)$. So, is it possible to construct a quantization of $\pi$?  

By the above arguments, $\mu_2$ exists (and 
a choice of $\mu_2$ is unique up to adding an arbitrary bivector). So there arises the question 
of existence of $\mu_3$ etc., i.e. the question whether there are other obstructions. 

The answer to this question is yes and no. Namely, if you don't pick $\mu_2$ carefully, 
you may be unable to find $\mu_3$, but you can always pick $\mu_2$ so that $\mu_3$ exists, and there 
is a similar situation in higher orders. This subtle fact is a consequence 
of the following deep theorem of Kontsevich:

\begin{theorem}\label{kon} 
\cite{K} Any Poisson structure $\pi$ on $A_0$ can be quantized. 
Moreover, there is a natural bijection between products $\mu$ 
up to an isomorphism and Poisson brackets $\pi_0+\hbar \pi_1+\hbar^2\pi_2+...$, 
such that the quasiclassical limit of $\mu$ is $\pi_0$.
\end{theorem}

{\bf Remark.} Note that, as was shown by O. Mathieu, 
a Poisson bracket on a general commutative $\Bbb C$-algebra may fail 
to admit a quantization. 

Let us consider the special case of symplectic manifolds, 
i.e. the case when $\pi$ is a nondegenerate bivector. 
In this case we can consider $\pi^{-1}=\omega$, which is a closed, nondegenerate 2-form 
(=symplectic structure) on $M$. In this case, Kontsevich's theorem  
is easier, and was proved by De Wilde - Lecomte, and later Deligne and Fedosov (see e.g. \cite{F}).
Moreover, in this case there is the following additional result, also due to Kontsevich, \cite{K}. 

\begin{theorem} If $M$ is symplectic and $A$ is a quantization of $M$, then 
the Hochschild cohomology $H^i(A[\hbar^{-1}])$ is isomorphic to $H^i(M,\Bbb C((\hbar)))$.  
\end{theorem} 

{\bf Remark.} Here the algebra $A[\hbar^{-1}]$ is regarded as a (topological) 
algebra over the field of Laurent series 
$\Bbb C((\hbar))$, so Hochschild cochains are, by definition, 
linear maps $A_0^{\otimes n}\to A_0((\hbar))$.

\begin{example}\label{quansym}
The algebra $B=A[\hbar^{-1}]$ provides an example of an algebra 
with possibly nontrivial $H^3(B)$, for which the universal deformation parametrized by $H=H^2(B)$ exists. 
Namely, this deformation is attached through the correspondence of Theorem \ref{kon}
(and inversion of $\hbar$)
to the Poisson bracket $\pi=(\omega+t_1\omega_1+...+t_r\omega_r)^{-1}$, where $\omega_1,...,\omega_r$ 
are closed 2-forms on $M$ which represent a basis of $H^2(M,\Bbb C)$, and $t_1,...,t_r$ are the 
coordinates on $H$ corresponding to this basis. 
\end{example}

\subsection{Examples} 

\begin{example}\label{refl} Let $V$ be a symplectic vector space over $\Bbb C$ with symplectic form $\omega$. 
Let ${\rm Weyl}(V)$ denote the Weyl algebra 
of $V$, which is the quotient of the free (=tensor) algebra on $V$ by 
the ideal generated by elements \linebreak $xy-yx-\omega(x,y)$. 

Let $G$ be a finite group acting symplectically on $V$. 
Then $G$ acts on ${\rm Weyl}(V)$, and one can form a semidirect product
algebra $A_0=G\ltimes {\rm Weyl}(V)$. Let us study deformations of $A_0$. 

We say that an element $g\in G$ is a symplectic reflection
in $V$ if \linebreak ${\rm rank}(g-1)|_V=2$. Let $S$ be the set of symplectic reflections
in $G$. 

\begin{proposition}\cite{AFLS}
$H^i(A_0)$ is the space of functions 
on the set of conjugacy classes of elements $g\in G$ such that 
${\rm rank}(g-1)|_V=i$. In particular, $H^i(A_0)=0$ if $i$ is odd, 
and $H^2(A_0)=\Bbb C[S]^G$. 
\end{proposition}

\begin{corollary}\label{sra}
There exists a universal deformation $A_u=\bold H_c(V,G)$ of $A_0$, which is parametrized 
by $c\in \Bbb C[S]^G$. 
\end{corollary}

The algebra $\bold H_c(V,G)$ is called the 
symplectic reflection algebra (see \cite{EG}).  
Such algebras were first considered by Drinfeld in 1986. 
If $V=\h\oplus \h^*$, where $\h$ is a representation of $G$, 
and the symplectic form on $G$ is the pairing between $\h$ and $\h^*$, 
then $\bold H_c(V,G)$ is called the rational Cherednik algebra.
We will later construct $\bold H_c(V,G)$ explicitly. 
\end{example}

\begin{example}\label{cher}
Let $X$ be a smooth affine algebraic variety over $\Bbb C$,
with an action of a finite group $G$. Let $D(X)$ be the algebra 
of algebraic differential operators on $X$. Let $A_0=G\ltimes D(X)$. 
Let us study deformations of $A_0$. 

For every $g\in G$, the fixed set $X^g$ of $g$ in $Y$ is a smooth affine variety, 
which consists of connected components $X^g_j$, possibly of different dimensions. 
Such a component is said to be a {\bf reflection hypersurface} if 
it has codimension $1$ in $X$. Let $S$ be the set of pairs 
$(g,Y)$, where $g\in G$, and $Y\subset X^g$ is a connected component which is a reflection hypersurface 
(i.e., has codimension 1). 
 
\begin{proposition} \cite{E}
One has $H^2(A_0)=(H^2(X,\Bbb C)\oplus \Bbb C[S])^G$. Moreover, 
there exists a universal deformation of $A_0$ parametrized by 
$H=H^2(A_0)$. 
\end{proposition}

This deformation $\bold H_c[X,G]$ is called the rational Cherednik algebra 
attached to $(X,G)$, and is described in \cite{E}.
If $X$ is a vector space $\h$ and $G$ acts linearly, 
then $\bold H_c[\h,G]=\bold H_c(\h\oplus \h^*,G)$ is the rational Cherednik 
algebra discussed above.  
\end{example}

\begin{example}\label{gen}
The following example from the paper \cite{DE} (conjecturally) 
generalizes examples \ref{quansym},\ref{refl}, and \ref{cher}.

Let $M$ be a symplectic $C^\infty$-manifold 
(or affine complex algebraic variety). Let $G$ be a finite group acting on 
$M$ by symplectic transformations, and $B$ be a quantization of 
$M$ which is equivariant under $G$ (such a quantization always exists). 
Let $A_0=G\ltimes B[\hbar^{-1}]$. Let us study deformations of $A_0$. 

The Hochschild cohomology of $A_0$ is given by the following theorem.
Let the fixed set $M^g$ be the union of connected components $M^g_i$, $i=1,...,N_g$. 

\begin{theorem}\label{coho} $H^*(A_0)$ equals, as a vector space, 
the orbifold cohomology of $M/G$ with coefficients in $\Bbb C((\hbar))$. 
Namely, 
$$
H^p(A_0)=
(\oplus_{g\in G}\oplus_{i=1}^{N_g}
H^{p-{\rm codim}M^g_i}(M^g_i))^G.
$$
(where the coefficients on the RHS are $\Bbb C((\hbar))$).  
\end{theorem}

{\bf Remark.} Let $S$ be the set of pairs $(g,Y)$, where $g\in G$, and 
$Y\subset M^g$ is a connected component of codimension $2$. 
Theorem \ref{coho} implies that 
$H^2(A_0)=(H^2(M)\oplus \Bbb C[S])^G$. 

Thus, we see that $H^3(A_0)$ does not always vanish.
Nevertheless, we make the following conjecture. 

\begin{conjecture}\label{conj}
The deformations of the algebra $A_0$ are
unobstructed. Thus there exists a universal 
deformation $H_c$ of this algebra parametrized by  
$c\in H^2(A_0)$. 
\end{conjecture}

Thus the conjecture implies that if $S\ne \emptyset$, then there exist
``interesting'' deformations of $A_0$, i.e., ones not coming from
$G$-invariant deformations of $B$. 

Let us give a few examples in which 
this conjecture is true. 

1. $H^3(A_0)=0$. This includes the following interesting case considered in \cite{EO}: 
$\Sigma$ is a smooth affine algebraic surface such that $H^1(\Sigma,\Bbb C)=0$,
and $M=\Sigma^n$, $G=S_n$. In this case there is 
one interesting deformation parameter 
corresponding to reflections in $S_n$. 

2. $G$ is trivial (Example \ref{quansym}).

3. $M=T^*Y$, where $Y$ is a smooth affine variety, 
and $G$ acts on $Y$ (Example \ref{cher}). 

4. If $M=V$ is a symplectic vector space and $G$
acts linearly (Example \ref{refl}).

5. Let $M=V/L$, where $V$ is a symplectic vector
space and $L$ a lattice in $V$ (i.e., $L$ is the abelian group
generated by a basis of $V$). Thus $M$ is an algebraic
torus with a symplectic form. We assume that the symplectic form
is integral and unimodular on $L$. Let $G\subset Sp(L)$ be a
finite subgroup; then $G$ acts naturally on $M$. 
In this case $H_c$ is an ``orbifold Hecke algebra'' 
defined in \cite{E} (it will be discussed below). 
\end{example}
    
\section{Algebras given by generators and relations}

\subsection{Giving formal deformations by generators and relations}

Another approach to exploring deformations of $A_0$ is 
defining deformations by generators and relations. 

Let us first consider the setting of formal deformations, which we 
have discussed in the previous section.  
Namely, let $A_0$ be an algebra over a field $k$, 
generated by $a_1,a_2,...$ with defining relations 
$R_j^0(a_1,a_2,...)=0$ (here $R_j^0$ 
are elements in the free $k$-algebra $F$ generated 
by $a_i$). Let us now define a formal deformation
of $A_0$ as the algebra over $K=k[[\hbar]]$ 
with the same generators and deformed relations $R_j=R_j^0+\hbar R_j^1+\hbar^2 R_j^2+...$. 
That is, $A$ is the quotient of the free algebra $F[[\hbar]]$ 
by the $\hbar$-adically closed ideal generated by the relations $R_j$. 

\begin{example} (The Weyl algebra.) Let $A_0=\Bbb C[x,y]$ be the algebra 
generated by $x,y$ with the defining relation $yx-xy=0$. 
We can then define $A$ by the same generators 
and the deformed relation $yx-xy=\hbar$ (the Heisenberg indeterminacy relation). 
Then $A$ is indeed a 1-parameter flat formal deformation of $A_0$, 
which provides a quantization of the standard Poisson bracket
$\lbrace{y,x\rbrace}=1$. 
\end{example}

So, is $A$ always a 
1-parameter flat formal deformation of $A_0$? In general the answer is {\bf no}:
the flatness property can fail. The following typical example of this
is obtained by adding just one relation to the relations above.  

\begin{example} Assume the algebra $A_0$ is defined by 
generators $x,y$ and defining relations
$$
yx-xy=0,\ x=0,
$$
and $A$ is defined by 
generators $x,y$ and relations 
$$
yx-xy=\hbar,\ x=0.
$$
Then $A$ is not topologically free, as it contains $\hbar$-torsion. Indeed,
$\hbar\cdot 1=yx-xy=0$ since $x=0$. On the other hand, 
$1\ne 0$, since the algebra $A_0=\Bbb C[y]$ is nonzero. 
\end{example}

In fact, it is easy to show that if we add any 
relation to $xy-yx=\hbar$, it will produce a non-flat deformation
(unless the algebra to be deformed is zero to begin with).
This shows that if one wants to secure flatness, one has to deform the relations 
in a very special way. In fact, it is usually rather difficult to 
do so, as well as to check that the resulting deformations are actually flat. 
Below I would like to show several situations when this task can be successfully completed. 
 
\subsection{Deformations of quadratic algebras}

The first situation is deformation theory of quadratic algebras.

Let $R$ be a finite dimensional semisimple algebra (say over $\Bbb C$). 
Let $A$ be a $\Bbb Z_+$-graded algebra, $A=\oplus_{i\ge 0}A[i]$, such that $A[0]=R$. 
For simplicity assume that the spaces $A[i]$ are finite dimensional for all $i$. 

\begin{definition} (i) The algebra $A$ is said to be quadratic if 
it is generated over $R$ by $A[1]$, and has 
defining relations in degree 2.

(ii) $A$ is Koszul if all elements of $Ext^i(R,R)$ 
(where $R$ is the augmentation module over $A$)
have grade degree precisely $i$. 
\end{definition}

{\bf Remarks.} 1. Thus, in a quadratic algebra, 
$A[2]=A[1]\otimes_{R} A[1]/E$, where $E$ is the 
subspace ($R$-subbimodule) of relations. 

2. It is easy to show that 
a Koszul algebra is quadratic, since 
the condition to be quadratic is just the Koszulity condition 
for $i=1,2$.

3. Many important algebras, e.g. the 
free algebra, the polynomial algebra and 
the exterior algebra are Koszul. 

Now let $A_0$ be a quadratic algebra, $A_0[0]=R$. 
Let $E_0$ be the space of relations for $A_0$.
Let $E\subset A_0[1]\otimes_{R} A_0[1][[\hbar]]$ 
be a topologically free (over $\Bbb C[[\hbar]]$) $R$-subbimodule which 
reduces to $E_0$ modulo $\hbar$ (``deformation of the relations''). 
Let $A$ be the ($\hbar$-adically complete) 
algebra generated over $R[[\hbar]]$ by $A[1]=A_0[1][[\hbar]]$ 
with the space of defining relations $E$. 
Thus $A$ is a $\Bbb Z_+$-graded algebra. 

Then we have the following fundamental result

\begin{theorem}\label{KDP} (Koszul deformation principle,\cite{D},\cite{BG},\cite{PP},\cite{BGS})
If $A_0$ is Koszul then $A$ is a topologically free $\Bbb C[[\hbar]]$ module
if and only if so is $A[3]$.  
\end{theorem} 

{\bf Remark.} Note that $A[i]$ for $i<3$ are obviously topologically free. 

\subsection{Symplectic reflection algebras.}
We will now demonstate by an example how the Koszul deformation principle works.

Let $V$ be a finite dimensional symplectic vector space 
over $\Bbb C$ with a symplectic form $\omega$, and $G$
be a finite group acting symplectically
on $V$. For simplicity let us assume that $(\wedge^2V)^G=\Bbb C\omega$. 

If $s\in G$ is a symplectic reflection, then let $\omega_s(x,y)$ be 
the form $\omega$ applied to the projections of $x,y$ 
to the image of $1-s$ along the kernel of $1-s$; thus $\omega_s$ is a skewsymmetric form of rank $2$ on $V$.   

Let $S\subset G$ be the set of symplectic reflections, and 
$c: S\to \Bbb C$ be a function which is invariant under the action of $G$. 
Let $t\in \Bbb C$. 

\begin{definition} The symplectic reflection algebra $H_{t,c}=H_{t,c}(V,G)$
is the quotient of the algebra $G\ltimes \bold T(V)$ 
by the ideal generated by the relation 
\begin{equation}\label{mainrel}
[x,y]=t\omega(x,y)-2\sum_{s\in S}c_s\omega_s(x,y)s.
\end{equation}
\end{definition}

The following theorem shows that the algebras 
$H_{t,c}(V,G)$ satisfy a flatness property, and moreover, they are 
the only ones satisfying this property within a certain natural class. 

\begin{theorem}\label{pbw1}
Let $\kappa: \wedge^2V\to \Bbb C[G]$ be a $G$-equivariant function ($G$ acts on the target by conjugation). 
Define the algebra $H_\kappa$ to be the quotient of the algebra 
$G\ltimes {\bold T}(V)$ by the relation 
$[x,y]=\kappa(x,y)$, $x,y\in V$. Put an increasing filtration on $H_\kappa$ 
by setting $\deg(V)=1$, $\deg(G)=0$, and define $\xi: G\ltimes SV\to {\rm gr}H_\kappa$ 
to be the natural surjective homomorphism. Then $\xi$ is an isomorphism if and only if 
$\kappa$ has the form 
$$
\kappa(x,y)=
t\omega(x,y)-2\sum_{s\in S}c_s\omega_s(x,y)s,
$$
for some $t\in \Bbb C$ and $G$-invariant function $c: S\to \Bbb C$. 
\end{theorem}

Before proving this theorem, let us point out a corollary. 
Denote by $\bold H_c(V,G)$ the algebra defined as $H_{t,c}(V,G)$, 
but with $t=1$ and $c$ being a formal parameter. 

\begin{corollary}
The algebra $\bold H_c(V,G)$ is a flat formal deformation 
of $G\ltimes {\rm Weyl}(V)$, parametrized by $\Bbb C[S]^G$. 
\end{corollary}

In fact, it turns out (see \cite{EG}) that $\bold H_c(V,G)$ 
is the universal deformation of $G\ltimes {\rm Weyl}(V)$, 
whose existence was proved in Example \ref{refl}. 

\begin{proof} (of Theorem \ref{pbw1})
Let $\kappa: \wedge^2V\to \Bbb C[G]$ be an equivariant map. 
We write $\kappa(x,y)=\sum_{g\in G}\kappa_g(x,y)g$, where $\kappa_g(x,y)\in \wedge^2V^*$.
To apply Theorem \ref{KDP}, let us homogenize our algebras. 
Namely, let $A_0=(G\ltimes SV)\otimes \Bbb C[u]$. Also 
let $\hbar$ be a formal parameter, and 
consider the deformation $A=H_{\hbar u^2\kappa}$ of 
$A_0$. That is, $A$ is the quotient of 
$G\ltimes {\bold T}(V)[u][[\hbar]]$ by the relations $[x,y]=\hbar u^2\kappa(x,y)$. 
This is a deformation of the type considered in Theorem \ref{KDP}, and 
it is easy to see that its flatness 
in $\hbar$ is equivalent to Theorem \ref{pbw1}. 
Also, the algebra $A_0$ is Koszul, because the polynomial algebra $SV$ is a Koszul algebra. 
Thus by Theorem \ref{KDP}, it suffices to show that 
$A$ is flat in degree 3. 

The flatness condition in degree 3 is ``the Jacobi identity''
$$
[\kappa(x,y),z]+[\kappa(y,z),x]+[\kappa(z,x),y]=0,
$$
which must be satisfied in $G\ltimes V$. In components, 
this equation transforms into the system of equations
$$
\kappa_g(x,y)(z-z^g)+\kappa_g(y,z)(x-x^g)+\kappa_g(z,x)(y-y^g)=0
$$
for every $g\in G$ (here $z^g$ denotes the result of the action of 
$g$ on $z$). 

This equation, in particular, implies that if $x,y,g$ are such
that $\kappa_g(x,y)\ne 0$ then for any $z\in V$ 
$z-z^g$ is a linear combination of $x-x^g$ and $y-y^g$. 
Thus $\kappa_g(x,y)$ is identically zero unless the rank of $(1-g)|_V$ is at most 2, i.e. 
$g=1$ or $g$ is a symplectic reflection. 

If $g=1$ then $\kappa_g(x,y)$ has to be $G$-invariant, 
so it must be of the form $t\omega(x,y)$, where $t\in \Bbb C$. 

If $g$ is a symplectic reflection, then
$\kappa_g(x,y)$ must be zero for any $x$ such that $x-x^g=0$. 
Indeed, if for such an $x$ there had existed $y$ with 
$\kappa_g(x,y)\ne 0$ then $z-z^g$ for any $z$ would be a multiple of $y-y^g$, which is impossible since 
$Im(1-g)|_V$ is 2-dimensional. This implies that $\kappa_g(x,y)=
-2c_g\omega_g(x,y)$, and $c_g$ 
must be invariant. 

Thus we have shown that if $A$ is flat (in degree 3) then $\kappa$ must have the form 
given in Theorem \ref{pbw1}. Conversely, it is easy to see that 
if $\kappa$ does have such form, then the Jacobi identity holds. 
So Theorem \ref{pbw1} is proved. 
\end{proof}

\subsection{Deformation of representations}

Another method of estabishing flatness of a deformation
$A$ of $A_0$ defined by generators and relations 
is showing that a given faithful representation $M_0$
of the algebra $A_0$ 
(for example, the regular representation) 
can be deformed (flatly) to a representation 
$M$ of $A$. In this case it follows automatically that 
$A$ is flat. Let us give two examples of situations where this method can be applied. 

\begin{example}\label{comman} (see \cite{E}). 
Let $X$ be a connected, simply connected complex manifold, 
and $G$ a discrete group of automorphisms of $X$. 
In this case the quotient $X/G$ is a complex orbifold. 
Let $X'\subset X$ be the set of points having trivial stabilizer (it is a nonempty open 
subset of $X$). Define the braid group $\widetilde G$ of the orbifold $X/G$ to be the 
fundamental group of the manifold $X'/G$ with some base point $x_0$. 
We have a surjective homomorphism $\phi: \widetilde G\to G$, 
which corresponds to gluing back the points which have a nontrivial stabilizer. 
Let $K$ be the kernel of this homomorphism. 

The kernel $K$ can be described by simple relations, corresponding to reflection hypersurfaces 
in $X$. Namely, given a reflection hypersurface $Y\subset X$, we have a conjugacy class $C_Y$ 
in $\widetilde G$ which corresponds to the loop in $X'/G$ which goes counterclockwise around $Y$. 
Let $T_Y$ be a representative of $C_Y$. Also, let $G_Y\subset G$ be the stabilizer of 
a generic point on $Y$; this is a cyclic group of some order $n_Y$. Then 
it follows from basic topology that 
the elements $T_Y^{n_Y}$ belong to $K$, and $K$ is the smallest normal subgroup 
of $\widetilde G$ containing all of them. In other words, the group $G$ is the quotient of 
the braid group $\widetilde G$ by the relations
\begin{equation}\label{Tn1}
T_Y^{n_Y}=1.
\end{equation}

Now let $A_0=\Bbb C[G]$, and let us define a deformation $A$ of $A_0$ 
to be the quotient of the group algebra of the braid group $\widetilde G$ 
by a deformation of relations (\ref{Tn1}). Namely, for every reflection hypersurface
$Y\subset X$ we introduce formal parameters $\tau_{Y,j}$, $j=1,...,n_Y$
(which are conjugation invariant), and replace relations (\ref{Tn1}) by the relations 
\begin{equation}\label{Tn1a}
\prod_{j=1}^{n_Y}(T_Y-e^{2\pi ij/n_Y}e^{\tau_{Y,j}})=0.
\end{equation}

The quotient $A$ of $\Bbb C[\widetilde G][[\tau]]$ by these relations 
is called the {\bf orbifold Hecke algebra} of $X/G$, and denoted by ${\mathcal H}_\tau(X,G)$.

\begin{theorem}\label{fla}(\cite{E})
If $H^2(X,\Bbb C)=0$ then ${\mathcal H}_\tau(X,G)$ is a flat deformation of $\Bbb C[G]$. 
\end{theorem}

{\bf Remark.} If $X$ is $\Bbb C^n$ and $G=G_0\cdot L$, where $L$ is a lattice of rank $2n$ and $G_0$ 
is a finite group acting on $L$ then 
${\mathcal H}_\tau(X,G)$ is, essentially, the algebra which was mentioned in Example \ref{gen}. 

To illustrate the relevance of the condition $H^2(X,\Bbb C)=0$, let us 
consider the special case when $G$ is the triangle group $F_{p,q,r}$, 
generated by $a,b,c$ with defining relations 
$$
a^p=1,\ b^q=1,\ c^r=1,\ abc=1,
$$
where $p,q,r>1$ are positive integers. 
The group $G$ is the group generated by rotations around the vertices 
of a triangle with angles $\pi/p,\pi/q,\pi/r$, by twice the angle at the vertex. 
Let $S=\frac{1}{p}+\frac{1}{q}+\frac{1}{r}$. The triangle lies on the sphere, Euclidean plane, or hyperbolic plane
$X$ when $S>1$, $S=1$, and $S<1$, respectively. The deformation ${\mathcal H}_\tau(X,G)$ is generated by 
$a,b,c$ with defining relations   
$$
\prod_{j=1}^p (a-\alpha_j)=0,\ \prod_{j=1}^q (b-\beta_j)=0,\ \prod_{j=1}^r (c-\gamma_j)=0,\ abc=1,
$$
where 
$$
\alpha_j=e^{2\pi ij/p}e^{\tau_{1j}},
\beta_j=e^{2\pi ij/q}e^{\tau_{2j}},
\gamma_j=e^{2\pi ij/r}e^{\tau_{3j}}.
$$
Theorem \ref{fla} says that the deformation is flat for the Euclidean and hyperbolic plane,
but says nothing about the sphere, i.e. the triples 
$(p,q,r)$ equal to $(2,2,n)$, $(2,3,3),(2,3,4),(2,3,5)$, 
 in which case the group $G$ is finite. 
And indeed, in this case ${\mathcal H_\tau}(X,G)$ is actually not flat!
To see this, note that in the sphere case ${\mathcal H}_\tau$, if it were flat, would have dimension $|G|$ 
(over $\Bbb C[[\tau]]$). So we may take the determinant of the relation $abc=1$ (using 
the fact that the eigenvalues of $a,b,c$ are $\alpha_j,\beta_j,\gamma_j$, with equal multiplicities). 
This yields a nontrivial relation on $\tau$:
$$
(\prod_{j=1}^p\alpha_j)^{|G|/p}
(\prod_{j=1}^q\beta_j)^{|G|/q}
(\prod_{j=1}^r\gamma_j)^{|G|/r}=1,
$$ 
which rules out flatness of ${\mathcal H}_\tau$. 
\end{example}

\begin{example}(\cite{ER}) Let $W$ be a Coxeter group of rank $r$
with generators $s_i$ and defining relations 
$$
s_i^2=1,\ (s_is_j)^{m_{ij}}=1\text{ for }m_{ij}<\infty,\ i,j=1,...,r,\ i\ne j,
$$
where $m_{ij}=m_{ji}$ are integers $\ge 2$ or $\infty$, defined for $i\ne j$.
Let $W_+$ be the subgroup of even elements of $W$. It is easy to see that 
$W_+$ is generated by the elements $a_{ij}:=s_is_j$, with defining relations 
$$
a_{ij}a_{ji}=1,\ a_{ij}a_{jk}a_{ki}=1,\ a_{ij}^{m_{ij}}=1.
$$ 
Define a deformation of $A_0=\Bbb C[W_+]$ as follows. 
Introduce invertible parameters $t_{ij,k}=t_{ji,-k}^{-1}$, $k\in \Bbb Z/m_{ij}\Bbb Z$
for $m_{ij}<\infty$. Let $R=\Bbb C[t_{ij,k}]$, and $A$ be the $R$-algebra 
generated by $a_{ij}$ with defining relations
$$
a_{ij}a_{ji}=1,\ a_{ij}a_{jk}a_{ki}=1,\ \prod_{k=1}^{m_{ij}}(a_{ij}-t_{ij,k})=0.
$$ 

For any $x\in W_+$, 
fix a reduced word $w(x)$ representing $x$. Let $T_{w(x)}$ be the element 
of $A$ corresponding to this word. 

\begin{theorem}\cite{ER}\label{er}
(i) The elements $T_{w(x)}$ for $x\in W_+$ span $A$ over $R$. 

(ii) These elements form a basis of $A$ over $R$ 
if and only if $W$ has no finite parabolic subgroups of rank $3$,
i.e. iff for each $i,j,l$, 
$$
\frac{1}{m_{ij}}+\frac{1}{m_{jl}}+\frac{1}{m_{li}}\le 1.
$$
\end{theorem}

\begin{corollary}
Let $\widehat A$ be the completion of $A$ with respect to the ideal generated by 
$t_{ij,k}-e^{2\pi k\sqrt{-1}/m_{ij}}$. Then $\widehat A$ is a flat deformation of $A_0$ 
iff $W$ has no finite parabolic subgroups of rank 3. 
\end{corollary}

{\bf Remark.} Note that triangle groups $F_{p,q,r}$ are groups $W_+$ for Coxeter groups 
of rank 3 (with $m_{12}=p$, $m_{23}=q$, $m_{31}=r$), 
so the ``only if'' part of Theorem \ref{er} (and the ``if'' part in rank 3) 
follow from Example \ref{comman}. 
\end{example}

In both of these examples, flatness is established 
by showing, using geometric methods (D-modules 
or constructible sheaves), that the regular representation 
of $A_0$ can be flatly deformed to a representation of the deformation.
Let us conclude by illustrating this in Example 
\ref{comman}, in the case when $X=E$ is a complex vector space, 
and $G$ is a finite group acting linearly on $E$. In this case, Theorem \ref{fla}  
was proved by Brou\'e, Malle, and Rouquier \cite{BMR},
following an idea of Cherednik. 
Let us sketch their proof. 

The main idea of the proof is to introduce Dunkl operators 
$D_a$, $a\in E$, which act on functions on $E$ (with poles 
on the reflection hyperplanes $Y$):
$$
D_a=\partial_a+\sum_Y\frac{\alpha_Y(a)}{\alpha_Y}(\sum_{g\in G_Y}c_{Y,g}g),
$$
where the summation is over all reflection hyperplanes $Y$, 
$\alpha_Y$ is the nonzero element of $E^*$ vanishing on $Y$, and 
$c_{Y,g}$ is a conjugation invariant function of $Y,g$. 

It can be shown that the Dunkl operators commute: $[D_a,D_b]=0$. 
This implies that the system of equations $D_a\psi=0$, $a\in E$, can be regarded as 
a local system with fiber $\Bbb CG$ on $(E\setminus \cup Y)/G$. 
The fundamental group of $(E\setminus \cup Y)/G$, by definition, is $\widetilde G$, so 
we may consider the corresponding monodromy representation of this group. 
If $c=0$, the monodromy representation is the standard homomorphism  
$\Bbb C \widetilde G\to \Bbb C G$. One may show that if $c\ne 0$, then the monodromy representation  
is a deformation of this standard homomorphism, which factors 
through the Hecke algebra ${\mathcal H}_\tau(E,G)$, for an appropriate linear change  
of variables $c\to \tau$. This implies the flatness of ${\mathcal H}_\tau(E,G)$.


\begin{thebibliography}{AFLS}

\bibitem[AFLS]{AFLS} J. Alev, M.A. Farinati, T. Lambre, and A.L. Solotar:
{\it Homologie des invariants d'une alg\`ebre de Weyl sous l'action d'un
groupe fini.} J. of Algebra {\bf 232} (2000), 564--577.

\bibitem[BG]{BG} A. Braverman, D. Gaitsgory,
{\it Poincar\'e-Birkhoff-Witt theorem for quadratic algebras of Koszul type},
J. Algebra {\bf 181} (1996), no.2, 315--328, {\tt{hep-th/9411113}}.

\bibitem[BGS]{BGS} A. Beilinson, V. Ginzburg, W. Soergel:
     {\it Koszul duality patterns in Representation Theory}, J. 
Amer. Math. Soc.
{\bf 9} (1996), 473--527.

\bibitem[BMR]{BMR} M. Brou\'e, G. Malle, and R. Rouquier, 
Complex reflection groups, braid groups, Hecke algebras.
J. Reine und Angew. Math. 500 (1998), 127-190.

\bibitem[DE]{DE} V. Dolgushev and P. Etingof, 
 Hochschild cohomology of quantized symplectic orbifolds and the Chen-Ruan cohomology, math.QA/0410562.

\bibitem[D]{D} V. Drinfeld: {\it
   On quadratic commutation relations in the quasiclassical case},
Mathematical physics, functional analysis (Russian),  25--34, 143, "Naukova Dumka", Kiev, 1986.
Selecta Math. Sovietica. {\bf 11} (1992), 317--326.

\bibitem[E]{E} P. Etingof, Cherednik and Hecke algebras of varieties with a finite
group action, math.QA/0406499.

\bibitem[EG]{EG} P. Etingof, V. Ginzburg, Symplectic reflection algebras, 
Calogero-Moser systems, and a deformed Harish-Chandra
isomorphism, Inventiones Math, vol. 147, (2002), p. 243-348.

\bibitem[EO]{EO} P. Etingof and A. Oblomkov,
Quantization, orbifold cohomology, and Cherednik algebras,
math.QA/0311005.

\bibitem[ER]{ER} P. Etingof, E. Rains, New deformations of group algebras of
Coxeter groups,
math.QA/0409261.

\bibitem[F]{F} B.V. Fedosov, A simple geometrical construction of deformation
quantization, J. Diff. Geom. {\bf 40} (1994) 213--238.

\bibitem[Ge]{Ge}
 Gerstenhaber, Murray On the deformation 
of rings and algebras.  Ann. of Math. (2)  79  1964 59--103.

\bibitem[HKR]{HKR} G. Hochschild, B. Kostant, and A. Rosenberg,
Differential forms on regular affine algebras, Trans. Amer. Math.
Soc. {\bf 102} (1962) 383-408.

\bibitem[K]{K} M. Kontsevich, Deformation quantization of Poisson
manifolds, Lett. Math. Phys. {\bf 66} (2003) 157-216;
q-alg/9709040.

\bibitem[PP]{PP} 
A. Polishchuk, L. Positselski, Quadratic algebras, Preprint, 1996;
to be published by the AMS, 2005. 

\end{thebibliography}
\end{document}